\documentclass[11pt]{amsart}
\usepackage{amsthm,amssymb,amsmath,calc,eucal,ifthen}

\newtheorem{theorem}{Theorem}[section]
\newtheorem{lemma}[theorem]{Lemma}
\newtheorem{proposition}[theorem]{Proposition}
\newtheorem{corollary}[theorem]{Corollary}

\newtheorem{remark}[theorem]{Remark}

\newtheoremstyle{rmdefinition}{}{}{\upshape}{}{\bfseries}{.}{
}{}

\theoremstyle{rmdefinition}

\newcommand{\derived}[2][1]{\ifthenelse{\equal{#1}{1}}{{#2}^{\prime}}{\ifthenelse{\equal{#1}{2}}{{#2}^{\prime\prime}}{\ifthenelse{\equal{#1}{3}}{{#2}^{\prime\prime\prime}}{{#2}^{(#1)}}}}}

\newcommand{\calA}{\mathcal{A}}

\newcommand{\calH}{\mathcal{H}}

\newcommand{\bbC}{\mathbb{C}}

\newcommand{\bbR}{\mathbb{R}}

\newcommand{\C}{{\mathbb C}}

\newcommand{\ol}[1]{\overline{#1}}
\newcommand{\interior}[1]{\overset{o}{#1}}
\newcommand{\rf}[1]{\cite{#1}}

\begin{document}

\bibliographystyle{plain2}

\title{Elliptic
operators on infinite graphs}

\author{J.~Dodziuk}
\thanks{This work is supported in part by a PSC-CUNY Research Grant.}

\address{ 
Ph.D.\ Program in Mathematics \\
Graduate Center (CUNY)\\
New York, NY 10016\\  
Email: jozek@derham.math.qc.edu}
\maketitle

\begin{abstract}
We present some 
applications of ideas from partial differential equations  and differential 
geometry to the study of difference 
equations on infinite graphs. All operators that we consider are 
examples of "elliptic operators" as defined by Y. Colin de Verdiere
\cite{CdV2}.
For such operators, we discuss analogs of inequalities of Cheeger and Harnack
and of the maximum principle (in both elliptic and 
parabolic versions), and apply them to study spectral theory, the ground 
state and the heat semigroup associated to these operators.
\end{abstract}

\section{Preliminaries}
We consider graphs (without loops or multiple connections) $G=(V,E)$
where $V$ is a set whose elements are called vertices and $E$, the set
of edges, is a subset of the set of two-element subsets of $V$. For an
edge $e=\{x,y\}\in E$, we will denote by $[x,y]$ the \emph{oriented}
edge from $x$ to $y$ and write $\overline{E}$ for the set of all
oriented edges. We also write $x \sim y$ if $\{x,y\}$ is an edge.
All graphs considered will be connected.

By a function on a graph we will mean a mapping $f:V\longrightarrow \C$.
By an operator on a graph, we shall always mean an operator acting on
functions and follow \cite{CdV2} in defining  the notion of
``self-adjoint, positive, elliptic operator.'' Observe first
that every operator $L$ is given by a matrix $(b_{x,y})$. We require our
operators to be local, i.e. $$
b_{x,y} =0 \qquad \text{if $\{x,y\}$ is not an edge and $x\neq y$
.}$$
Thus $$
Lf(x)=b_{x,x}f(x)+\sum_{x\sim y}b_{x,y} f(y) . 
$$
The constant functions are annihilated by $L$ if and only if 
$\sum_{y\sim x}b_{x,y} = - b_{x,x}$ for every $x\in V$. 
Every local
operator L can be rewritten in the form 
\begin{equation}\label{local}
Lf(x)=W(x)f(x)+\sum_{y\sim x}a_{x,y}(f(x)-f(y))
\end{equation}
where $W(x)=b_{x,x} + \sum_{y\sim x}b_{x,y}$ and 
$a_{x,y}=-b_{x,y}$.
We will often write $L=A+W$, where $A$, given by the sum in the
formula above, annihilates constant functions and $W$ denotes the operator of multiplication by the
function $W(x)$. 

Let $\ell^2(V)$ be the space of complex-valued functions $f$ satisfying 
$$\sum_{x\in V} |f(x)|^2 < \infty$$ equipped with the standard hermitian
inner product
$$
(f,g)=\sum_{x\in V} f(x)\ol{g(x)}.$$
We denote by $C_0(V)$ the space of all functions on $V$
with finite support. In order that the operator $L$ be symmetric on
$C_0(V)$, i.e.\ $(Lf,g) = (f,Lg)$ it is necessary and sufficient that 
$a_{x,y} = \overline{a_{y,x}}$ and $W(x)+\sum_{y\sim x}a_{x,y}\in \bbR$. 
We want to think of the operator in
(\ref{local}) as a ``Laplacian'' plus a potential. Thus, we impose an
additional condition on $A$ that will make it positive on $C_0(V)$.
Namely, we require that $a_{x,y}$ be real and positive for every 
edge $\{x,y\}$. We will
refer to such operators as \emph{elliptic}, \emph{positive} and \emph{symmetric}.
A very important example is the combinatorial Laplacian $A=\Delta$ given
by choosing $a_{x,y}=1$ for every edge,
$$ \Delta f(x) = \sum_{x\sim y}(f(x)-f(y))=
m(x)f(x) -\sum_{x\sim y}f(y),$$
where $m(x)$ is the valence of the vertex $x\in V$ i.e.\ the number of
edges emanating from $x$.

The following lemma sheds some light on the structure of a positive,
symmetric operator. First, we need a definition. Let $C(\overline{E})$ 
denote the
space of functions $\phi$ on \emph{oriented} edges satisfying
$\phi([x,y])=-\phi([y,x])$ for every edge $\{x,y\}$ and let $$
\ell^2(\overline{E}) = \{ \phi \in C(\ol{E}) \mid
\sum_{\{x,y\}\in E} |\phi([x,y])|^2 < \infty\}.$$ We equip 
$\ell^2(\overline{E})$ with the natural inner product
$$ <\phi,\psi> = \sum_{\{ x,y\}\in E}\phi([x,y])\ol{\psi([x,y])}.$$
In addition, given a positive, symmetric operator $A$ as above, define the (possibly unbounded) operator
$d_A$ from $\ell^2(V)$ to $\ell^2(\ol{E})$ by $$
d_A f([x,y]) =  \sqrt{a_{x,y}} (f(x) - f(y)).$$

\begin{lemma}{\label{divergence}}
Suppose $f$ and $g$ are two functions on the graph and one of them has 
finite support. Then $$
(Af,g) = <d_A f, d_A g>.$$
In particular, if $f$ has finite support, $(Af,f) \geq 0$ with equality
if and only if $f\equiv 0$.
\end{lemma}

\begin{proof}
The proof is a simple calculation. 
\begin{eqnarray*}
(Af,g)&=&\sum_x \left ( \sum_{y\sim
x} a_{x,y} (f(x) - f(y)) \right )\ol{g(x)}\\
& =&
\sum_{\{x,y\}\in E} a_{x,y} (f(x)-f(y))\ol{(g(x)-g(y))} = <d_A f,d_A g>
\end{eqnarray*}
To justify it note that an edge $\{z,w\}$ contributes to the first sum
twice. The contribution is 
\begin{eqnarray*}
a_{z,w} (f(z)-f(w))\ol{g(z)} &+& a_{w,z} (f(w)-f(z))\ol{g(w)} = \\[5pt] &&
a_{z,w} (f(z)-f(w))\ol{(g(z)-g(w))}
\end{eqnarray*}
since $a_{z,w}$ is symmetric. This
proves that the two sums are equal. The statement about strict
positivity of $(Af,f)$ follows trivially.
\end{proof}

We wish to consider $L=A+W$ as an unbounded operator on $\ell^2(V)$ and
to study its spectrum. In order to obtain a reasonable setup we will
require that the potential $W$ be bounded from below by a constant,
$W(x) \geq c$ for all $x \in V$. By the lemma above, $L$ is
semi-bounded, i.e.\ $(Lf,f) \geq c (f,f)$ for every $f\in C_0(V)$.
By Theorem X.23 of \cite{simon-reed2}, $L$ then has a distinguished 
self-adjoint extension, the Friedrichs extension, $\hat{L}$ such that 
$\lambda_0(\hat{L})$, 
the bottom of the spectrum of $\hat{L}$, has a variational
characterization 
\begin{equation}\label{var-princ}
\lambda_0(\hat{L}) = \inf_{f\in C_0(V)\setminus \{0\}}
\frac{(Lf,f)}{(f,f)}.
\end{equation}
We will abuse the notation and write $\lambda_0(L)$ for  
$\lambda_0(\hat{L})$.

In general, without any further restrictions, the operator $L$ with domain $C_0(V)$ may have many self-adjoint extensions. The theorem below gives
 conditions under which $L$ is essentially self-adjoint, i.e.\ has a unique self-adjoint extension, cf. \cite{simon-reed2}, Theorem X.28.

\begin{theorem}\label{a+w-sa}
Suppose that $A$ is a positive, symmetric and bounded as an operator on $\ell^2(V)$. Let $W$ be bounded from below by a constant. Then $L=A+W$ is essentially self-adjoint on $C_0(V)$. 
\end{theorem}
\begin{proof}
Choose a positive constant $\kappa$ so that $W+\kappa\geq 1$. By Theorem X.26 of \cite{simon-reed2}, it suffices to show that \begin{equation}
(A+W+\kappa)^*f=0 \label{adj-eq}
\end{equation}
implies that $f=0$. Taking the inner product of the equation above with the function $\delta_x$ 
($\delta_x(y) = 1$ if $x=y$ and $0$ otherwise), 
using the definition of the adjoint and Lemma \ref{divergence}, we see that (\ref{adj-eq}) is equivalent to $$
(A+W+\kappa)f=0, \qquad\qquad f\in \ell^2(V)$$
where $(A+W+\kappa)f$ is computed pointwise as in (\ref{local}) with $W$ replaced by $W+\kappa$.
Since $A$ is bounded and $C_0(V)$ is dense in $\ell^2(V)$, $(Af,f)\geq 0$ by Lemma \ref{divergence}. Therefore,
$0=(Af,f)+((W+\kappa)f,f) \geq (f,f)$. It follows that $f=0$ which proves the theorem.
\end{proof}
\begin{remark}\label{bounded}
Observe that the condition that $A$ be bounded holds if $a=\sup{a_{x,y}}<\infty$ and $M=\sup m(x) < \infty$. In fact,
in this case $\parallel A \parallel \leq 2 aM$.
\end{remark}
\begin{remark}
We view the Theorem \ref{a+w-sa} as an analog of Theorem X.28 of \cite{simon-reed2} which applies to a differential  operator $-\Delta + V$ on $\bbR^n$. Clearly, $\Delta$ is unbounded but the unboundedness is an infinitesimal effect that does not occur for difference operators on graphs. We view the boundedness of $A$ or the
condition $a<\infty$ as a partial replacement of uniform ellipticity, (see Corollary \ref{growth-gr-state} below for a proper analog of
uniform ellipticity). Similarly, $M< \infty$ is a bounded geometry condition.

\end{remark}

We now state two local results. Their continuous analogs - the maximum principle and Harnack's inequality -
are discussed at great length in \cite{protter-weinberger}. Let $V_1 \subset V$ be a set of vertices and let $G_1$ be the full subgraph of $G$ generated by $V_1$ (i.e.\ the set of edges of $G_1$ consists of all edges $\{x,y\}$ of G such that $x,y\in V_1$).
Let $\interior{V_1} = \{x\in V_1\mid y \sim x \quad \text{implies}\quad y\in V_1\}$ and $\partial V_1 = V_1\setminus
\interior{V_1}$. We say that $\interior{V_1}$ is connected if every two of its vertices $x$, $y$ can be connected by a
path of edges $[x_0,x_1],[x_1,x_2],\: \ldots\: ,[x_{n-1},x_n]$, $x_0=x, x_n=y$ with $x_i \in \interior{V_1}$ for
$i=0,1,\ldots ,n$.

\begin{lemma}\label{max-pr}
Let $L=A+W$ where $A$ is positive, symmetric and $W$ is nonnegative. 
Suppose $V_1 \subset V$ is a subset with $\interior{V_1}$ connected. Let $f$ be a function on $V_1$
such that $$Lf(x)=Af(x) +W(x)f(x) \geq 0\quad \text{for}\quad x\in \interior{V_1}.$$
If $f$ has a minimum at $x_0\in \interior{V_1}$ and $f(x_0)\leq 0$ then $f$ is constant on $V_1$.
\end{lemma}

\begin{proof}
Suppose $x_0\in \interior{V_1}$ is a minimum and $f(x_0)\leq 0$. Then $$
0\leq \sum_{y\sim x_0} a_{x_0,y}(f(x_0) -f(y)) +W(x_0)f(x_0) \leq 0
$$
since $A$ is positive, $x_0$ is a minimum, and $W(x_0)f(x_0) \leq 0$. It follows that all terms in
the sum above are equal to zero, i.e.\ $f(y) = f(x_0)$ for every $y\sim x_0$. By connectedness, $f$ is
constant.
\end{proof}

\begin{lemma}\label{harnack}
Suppose $A$ and $W$ satisfy the assumptions of Lemma \ref{max-pr}. Let $V_1\subset V$, $x\sim y$,
$x,y\in\interior{V_1}$. If $$
Lf=Af + Wf \geq 0\quad \text{and}\quad f>0 \quad \text{on}\quad V_1$$
then $$
\frac{a_{x,y}}{\left ( W(x) + \sum_{z\sim x} a_{x,z} \right )} \leq \frac{f(x)}{f(y)}\leq
\frac{\left ( W(y) + \sum_{z\sim y} a_{y,z} \right )}{a_{x,y}}.
$$
\end{lemma}
\begin{proof}
By symmetry, it suffices to prove one of the two inequalities above. We have
$$
(A+W)f(x) = \sum_{z\sim x} a_{x,z} (f (x) -f (z)) +W(x)f (x) \geq 0.$$
Therefore, $$
\left (\sum_{z\sim x} a_{x,z} \right )f(x) + W(x)f(x) \geq  \sum_{z\sim x} a_{x,z} f(z) 
\geq a_{x,y} f(y).$$
This, of course, is equivalent to the lower bound on $f(x)/f(y)$ in the statement of the lemma.
\end{proof}
We refer to Lemma \ref{max-pr} as the maximum principle and to Lemma \ref{harnack} as the Harnack inequality.
The significance of the Harnack inequality is that it gives a bound of the ratio $f(x)/f(y)$ in terms of 
the coefficients of the operator but \emph{independent} of the function $f$.

\section{Existence of ground state}
In this section we prove, for an operator $L=A+W$ with positive, symmetric A and the potential W bounded from below by a constant, the existence of a ground state, i.e.\ a positive solution of the equation $$L\phi =\lambda_0(L) \phi,$$
cf.\ \cite{pinsky} for an extensive discussion in the continuous setting.
We assume that the underlying graph $G$ is connected and fix a vertex $x_0$ as an ``origin''. Consider the exhaustion $\{G_n\}_{n=1}^\infty$ of $G$
where, for every $n$, $G_n$ is the full subgraph with the vertex set $V_n=\{x\in V \mid d(x_0,x) \leq n\}$. 
Here, $d(x,y)$ denotes the combinatorial distance between $x,y\in V$, i.e.\ the length of the shortest path of
edges connecting $x$ with $y$.
Clearly,
$\interior{V_n}$ is connected for every $n\geq 1$. We will construct a ground state $\phi$ by solving certain
``boundary value problems'' on $G_n$ and taking a limit of the solutions. In order to get started we need to review
these boundary value problems. Thus, let $U$ be a finite subset of $V$ such that the full subgraph generated by $U$ has connected interior. Let $C_0(U)$ be the space of functions on $U$ that vanish on $\partial U$.  Extending functions in
$C_0(U)$ by zero embeds $C_0(U)$ isometrically in $C_0(V)$. 
We define, for $f\in C_0(U)$, $L_Uf\in C_0(U)$ by 
\begin{equation*}
L_Uf(x)=\begin{cases}
            W(x)f(x)+\sum_{x\sim y}a_{x,y}(f(x)-f(y)) & \text{if $x\in \interior{U}$},\\
        0 & \text{if $x\in \partial U$}.
      \end{cases}
\end{equation*}
We can define $A_Uf\in C_0(U)$ for $f\in C_0(U)$ analogously. 
The calculation in the proof of 
Lemma \ref{divergence} shows that $A_U$ and $L_U$ are 
symmetric operators on $C_0(V)$ and that $A_U$ is strictly positive. It follows that $\lambda_0(L_U)$, the smallest eigenvalue of $L_U$ on $C_0(U)$, has variational characterization
\begin{equation}\label{var}
\lambda_0(L_U) = \inf_{f\in C_0(U)\setminus \{0\}}
\frac{(L_Uf,f)}{(f,f)}  =
\inf_{f\in C_0(U)\setminus \{0\}} \frac{(Lf,f)}{(f,f)} 
\end{equation}
where in the last expression above we identify $f$ with its extension 
by zero outside $U$.
\begin{proposition}
The eigenspace of $\lambda_0(L_U)$ is one-dimensional and every eigenfunction $\psi$ belonging to $\lambda_0(L_U)$ has constant sign in the interior of $U$.
\end{proposition}
\begin{proof}
It is enough to consider real-valued functions. Replacing $W$ by $W+c$ with a suitably large $c$, we can assume that $W$ is nonnegative. Since $$
(L_Uf,f) = \sum_{x\sim y, \: x\in  U,\:y\in\interior{U}} a_{x,y} (f(x)-f(y))^2 + \sum_{x \in \interior{U_0}}W(x)f(x)^2
$$
replacing $f$ by $|f|$ decreases the Rayleigh-Ritz quotient in (\ref{var}). Therefore, it follows that if $\psi$ is an
eigenfunction belonging to $\lambda_0(L_U)$ then $|\psi|$ is one as well. Thus we can assume that there exists a nonnegative eigenfunction $\psi$. Since the Raylegh-Ritz quotient is nonnegative, $\lambda_0(L_U) \geq 0$. The maximum principle in Lemma \ref{max-pr} implies that $\psi$ is strictly positive in $\interior{U}$. Finally, if the eigenspace of $\lambda_0(L_U)$ had two or more dimensions, there would exist another eigenfunction $\phi$ orthogonal to $\psi$. Therefore $\phi$ would have to change sign and be negative at an interior point, but this is impossible by the maximum principle.
\end{proof}
We are now ready to prove 
\begin{theorem}\label{gr-state}
Consider an operator $L=A+W$ on a connected graph $G$ with positive, symmetric A and the potential $W$ bounded below by a constant.
There exists a ground state $\phi$ for $L$ i.e.\ a function $\phi > 0$ on $V$ such that
$$L\phi = \lambda_0\phi$$
where $\lambda_0=\lambda_0(L)$ is the bottom of the spectrum of (the Friedrichs extension of) $L$ on G. 
\end{theorem}
\begin{proof}
The proof for the case of the combinatorial Laplacian was given in \cite{dod-mat3}. We follow the same line of
argument here but remark that exhaustion argument of this kind is applied very often in studying partial differential equations on noncompact domains or domains with non-smooth boundaries as, for example, in \cite{pinsky}, Chapter 4. Note first that by adding a suitable constant to the potential $W$ we can assume without any loss of
generality that $W>0$. We use the exhaustion of $G$ by finite subgraphs $G_n$ described above. Let $\lambda_n=\lambda_0(L_{G_n})$ and let $\phi_n$ be the corresponding positive eigenfunction of $L$ on $C_0(V_n)$ normalized so that $\phi_n(x_0)=1$. By the variational characterization of eigenvalues and of the bottom of the spectrum
(\ref{var-princ}), (\ref{var}) we have $\lambda_n\searrow \lambda_0$. Fix a point $y\in V$. Then, there exists
$k=k(y)$ such that $y\in \interior{V_n}$ for all $n>k$. Choose a path of length $d(x_0,y)$ that connects $x_0$ and $y$.
Using the normalization $\phi_n(x_0)=1$ and applying the local Harnack inequality in Lemma \ref{harnack} to successive edges of the path, we see that the sequence $\phi_n(y)$ is bounded above and below by positive constants that are independent of $n$. Using the diagonal process, we choose a subsequence $(n_k)_{k=1}^\infty$ such that
the sequence $(\phi_{n_k}(y))_{k=1}^\infty$ converges to the limit $\phi(y)$ of every vertex $y\in V$ and $\phi(y)> 0$. Since $L\phi$ is given by
the formula (\ref{local}) and $\lambda_n\searrow \lambda_0$ we see that $\phi$ is a positive solution of 
$L\phi=\lambda_0 \phi$ as required.
\end{proof}
We now need the following lemma to control the behavior at infinity of a ground state under certain additional assumptions.
\begin{corollary} \label{growth-gr-state}Assume that $A$ is symmetric
and positive, 
that the graph $G$ has bounded valence $\sup_{x\in V} m(x) = M <\infty$ and that the operator $A$ is uniformly elliptic in the sense that there exist constants 
$\gamma, \Gamma > 0$ so that $\gamma \leq a_{x,y} \leq \Gamma$ for every edge $\{x,y\}$. Suppose a function $f$ on $V$
satisfies $Af\geq 0$, $f>0$. Then, for every $x,y\in V$, $$
\left (\frac{M\Gamma}{\gamma}\right )^{-d(x,y)} \leq \frac{f(x)}{f(y)}  \leq \left (\frac{M\Gamma}{\gamma}\right )^{d(x,y)}.$$
\end{corollary}
\begin{proof}
By Lemma \ref{harnack}, $\gamma/M\Gamma \leq f(z)/f(w) \leq M\Gamma/\gamma$ if $z\sim w$. We connect $x$ with $y$ by a path of
edges of length $d(x,y)$ and apply  these inequalities for every edge along the path. The corollary follows.
\end{proof}
Observe that this is entirely analogous to Theorem 21 in \cite{protter-weinberger}.
\section{Cheeger's inequality}
In this section, we assume that $L=A$ and give a lower bound for the bottom of the spectrum of $A$ on $G$. This bound
originated in Riemannian geometry \cite{cheeger} and has been studied a great deal for the combinatorial Laplacian on graphs \rf{lubotzky}, \rf{d2}, \rf{dod-ken}. 

As before, let $A$ be a positive, symmetric elliptic operator on an infinite graph $G$ and let $U\subset V$ be a finite subset. We define
\begin{equation} \label{iso-U}
h_A(U) = \frac{\sum_{x\in\interior{U},\,y\in \partial U,\,x\sim y} \sqrt{a_{x,y}}}{\# (U)},
\end{equation}
 and 
\begin{equation}\label{cheeger}
\beta (G,A) = \inf_U h_A(U)
\end{equation}
where $\# U$ denotes the number of vertices of $U$. 
\begin{theorem}\label{cheeger-inq} Suppose $\sup_{x\in V} m(x) = M < \infty$
The lower bound of the spectrum of $A$ on $G$ satisfies $$
\lambda_0(A) \geq  \frac{\beta(G,A)^2}{2M}.$$
\end{theorem}
\begin{proof}
We follow the proof of Theorem 2.3 of \rf{d2}. Let $(G_n)_{n=1}^{\infty}$ be the exhaustion of $G$ used in
the proof of Theorem \ref{gr-state}. Since $\lambda_n\searrow \lambda_0$ it will suffice to show that
$\lambda_n \geq  \beta(G,A)^2/{2M}$ independently of $n$. We will fix $n$, set $U=V_n$ and let $\phi$ be positive
eigenfunction of $A_U$. Observe that by Lemma \ref{divergence} 
and (\ref{var})
\begin{equation}\label{lambda}
\lambda_n = \lambda_0(A_U) = \frac{<d_A \phi, d_A\phi>}{(\phi,\phi)}
\end{equation}
if we extend $\phi$ by zero outside $U$. Consider the expression $$
\calA = \sum_{\{x,y\}\in E} \sqrt{a_{x,y}}|\phi^2(x)-\phi^2(y)|.$$
By Cauchy-Schwartz inequality we have 
\begin{eqnarray*}
\calA &=& \sum_{\{x,y\}}\sqrt{a_{x,y}} |\phi(x)-\phi(y)|\,|\phi(x) + \phi (y)|\\
      &\leq & \left ( \sum_{\{x,y\}} |\phi(x)+\phi(y)|^2\right )^{1/2} \,\left (\sum_{\{x,y\}} a_{x,y} |\phi(x)-\phi(y)|^2\right )^{1/2}\\
      &\leq& \sqrt{2}\left ( \sum_{\{x,y\}} (\phi^2(x)+\phi^2(y))\right )^{1/2} \, (d_A\phi ,d_A\phi )^{1/2}.
\end{eqnarray*}
In $\sum_{\{x,y\}}(\phi^2(x) + \phi^2(y))$, every vertex contributes as many times as the number of edges emanating from it.
Hence we get the following upper bound on $\calA$.
\begin{equation}\label{a-upper}
\calA \leq \sqrt{2M}\,(\phi ,\phi)^{1/2}\, (d_A\phi, d_A\phi )^{1/2}.
\end{equation}
On the other hand we can estimate $\calA$ from below in terms of $(\phi,\phi)$ as follows. Let $0=\nu_0<\nu_1<\nu_2< \ldots<\nu_N$ be the sequence of all values of $\phi^2$. Note that, since $A\phi (x) = \lambda_0(U) \phi (x)$ at
every interior vertex $x$ and since $\lambda_0(U) > 0$ by (\ref{lambda}), every interior vertex $x$ will have a neighbor $y$ such that $\phi(x) > \phi (y)$. Define a set of vertices $U_i$, $i=1,2,\ldots ,N$  as follows.
A vertex $x\in U$ belongs to $U_i$ if and only if $\phi^2(x) \geq \nu_i$ and let $F_i$ be the full graph generated by
the set $U_i$. Now $$
\calA = \sum_{i=1}^N \sum_{\phi^2(x)=\nu_i} \sum_{y\sim x\,\phi^2(y)<\nu_i} \sqrt{a_{x,y}}(\phi^2(x)-\phi^2(y)).$$
If $\phi^2=\nu_i$ and $\phi^2(y)=\nu_{i-k}$ for some $k\in \{1,2,\ldots,i\}$, then on the one hand,$\phi^2(x)-\phi^2(y)=(\nu_i-\nu_{i-1}) + (\nu_{i-1} - \nu_{i-2} + \ldots (\nu_{i-k+1}-\nu_{i-k})$ and, on the other hand, $x \in \partial U_i \cap \partial U_{i-1} \cap \ldots \cap \partial U_{i-k+1}$. It follows that
$$
\calA\geq \sum_{i=1}^N (\nu_i-\nu_{i-1})\sum_{y\sim x,\,y\in\partial U_i} \sqrt{a_{x,y}}.$$
Applying (\ref{cheeger}) we obtain $$
\calA \geq h_A(U) \sum_{i=1}^N \#U_i(\nu_i - \nu_{i-1})\geq \beta \sum_{i=1}^N \#U_i(\nu_i - \nu_{i-1})$$
with $\beta=\beta(G,A)$.
``Summation by parts'' now yields$$
\calA \geq \beta\left ( \nu_N \#U_n + \sum_{i=1}^{N-1} \nu_i(\#U_i - \# U_{i+1})\right ) .$$
Observe that $\#U_n$ is the cardinality of the set where $\phi^2=\nu_N$ while $\#U_i-\#U_{i+1}$ is the number of points where $\phi^2=\nu_i$. It follows that $$
\calA \geq \beta (\phi,\phi).$$
This inequality combined with (\ref{lambda}) and (\ref{a-upper}) gives the desired lower bound.
\end{proof}

We remark that one can also bound $\lambda_0(A)$ from above by a related isoperimetric constant. Namely, let
$\chi_U$ be the characteristic function of a finite set of vertices $U\subset V$. Then $$
\lambda_0(A) \leq \frac{<d_A\chi_U,d_A\chi_U>}{(\chi_U,\chi_U)} = \frac{\sum_{x\sim y,x\in U, y\not\in U}a_{x,y}}{\#U}$$
It follows that $$\lambda_0(A) \leq \beta_1(G,A) = \inf \frac{\sum_{x\sim y,x\in U, y\not\in U}a_{x,y}}{\#U}$$
where the infimum is taken over all finite subsets $U$ of $V$.

Note that for the combinatorial Laplacian $\Delta$, $a_{x,y}\equiv 1$. Thus $\beta(G,A)=\beta_1(G,A)$.
In particular, for graphs of bounded valence $\lambda_0(\Delta) =0$ if and only if $\beta(G,\Delta) =0$ which is analogous to a result of Buser \cite{buser-uppr} in the Riemannian setting and is very useful in connection with
various characterizations of amenability of discrete, finitely generated groups \cite{brooks2}.

\section{The heat equation}
In this section we make several standing assumptions. Namely, we assume that the graph G has bounded valence
$\sup_{x\in V} m(x) =M < \infty$; that the potential $W\equiv 0$ i.e.\ $L=A$; and that $a=\sup_{\{x,y\}\in E} a_{x,y} <\infty$. We shall study the parabolic initial value problem 
\begin{equation}\label{initial}
\begin{split}
Au + \frac{\partial u}{\partial t}  &= 0\\
u(x,0) &= u_0(x)
\end{split}
\end{equation}
 and the associated heat semigroup
using the method of \cite{d3} applied previously to the combinatorial Laplacian in
\cite{dod-mat3}. Here $u(x,t)$ is a function of $x\in V$ and
$t>0$, while $u_0$ is a given function on $G$. The first equation above will be referred to
as the heat equation.

We are going to use the following version of the maximum principle, see \cite{protter-weinberger},
Chapter 4 for an analog in the continuous setting.
\begin{lemma}\label{max-par}
Suppose $u(x,t)$ satisfies the inequality $Au + \frac{\partial
u}{\partial t} < 0$ on $\overset{o}U \times [0,T]$ for a finite
subset $U$ of $V$. Then the
maximum of $u$ on $U\times [0,T]$ is attained on the set
$U\times \{0\} \cup \partial U\times [0,T]$.
\end{lemma}
\begin{proof}
Suppose $(x_0,t_0)\in \overset{o}{U} \times (0,T]$ is a maximum.
It follows that $\frac{\partial u}{\partial t}(x_0,t_0)$ is nonnegative so
that $Au(x_0,t_0) < 0$. On the other hand, (\ref{local}) and
positivity of $A$ imply that $Au(x_0,t_0) \geq 0$. The
contradiction proves the lemma.
\end{proof}
We use the lemma above to prove the uniqueness of bounded solutions of (\ref{initial}).
\begin{theorem}\label{unique}
Let $u(x,t)$ be a bounded solution of (\ref{initial}) with
the initial condition
$|u_0(x)|\leq N_0$. Then $u$ is determined uniquely by $u_0$
and $$ |u(x,t)| \leq N_0$$
for all $(x,t)$. Moreover, if a bounded initial condition $u_0$ is given,
then a bounded solution $u(x,t)$ of (\ref{initial}) exists.
\end{theorem}
\begin{proof}
Suppose that $u(x,t)$ is a bounded solution. Let $N_1= \sup
|u(x,t)|$. Fix $x_0 \in V$ and define $r(x) = d(x,x_0)$. By our
assumption on the valence and (\ref{local})
\begin{equation} \label{Ar}
|Ar| \leq aM .
\end{equation}
Consider an auxiliary function $$ v(x,t) = u(x,t) - N_0 -
\frac{N_1}{R} \left(r(x) + a(M+1)t\right ),$$ where $R$ is a large parameter. Let
$U = B(x_0,R)$ be the set of vertices of $V$ at distance at most
$R$ from $x_0$. The function $v(x,t)$ is nonpositive on the set
$U\times \{0\} \cup \partial U\times [0,T]$ and satisfies $(A
+ \frac{\partial}{\partial t}) v < 0$ on
$\overset{o}U \times[0,T]$ because of (\ref{Ar}). Lemma
\ref{max-par} implies therefore that $v(x,t) \leq 0$ so that $$
u(x,t) \leq N_0 + \frac{N_1}{R}\left (r(x) + a(M+1)t\right )$$ on $B(x_0,R)\times
[0,T]$. Keeping $(x,t)$ fixed and letting $R$ increase without
bounds, we see that $ u(x,t)\leq N_0$. Applying the same argument
to $-u$ yields $|u(x,t)| \leq N_0$. Since $T > 0$ and $x$ were
arbitrary, this last inequality holds for all $x\in V $ and $t\geq
0$. Uniqueness follows by considering the difference of two
solutions. We postpone the proof of existence of the solution.
\end{proof}

Recall that under our assumption $A$ is a bounded operator on $\ell^2(V)$.
Therefore, we can define for $t\geq 0$ \begin{equation}
\label{semi}
P_t = e^{-tA} = \sum_{k=0}^\infty (-1)^k \frac{t^kA^k}{k!}.
\end{equation}
Obviously, $u(x,t)=\left ( P_t u_0\right )(x)$ is a solution of (\ref{initial}) whenever $u_0$ is
in $\ell^2(V)$. Since $\parallel P_t\parallel \leq 1$ we see that 
for every $x \in V$ and $t\geq 0$ $$
|u(x,t)| \leq \parallel u(\cdot , t) \parallel \leq \parallel u_0 \parallel$$
so that $u(x,t)$ is a bounded solution and we get uniqueness.
We would like to extend the semigroup $P_t$ to a larger class of 
functions.

We define $p_t(x,y)$ to be matrix coefficients of the operator $P_t$, i.e. $$
p_t(x,y)=(P_t\delta_x,\delta_y)$$ 
where $\delta_x$ is the characteristic function of the set $\{x\}$.
Similarly, let $A(x,y)=(A\delta_x,\delta_y)$.
Since $A$ is self-adjoint both of these matrices are symmetric. Writing $u_0 = \sum_y u_0(y)\delta_y$
and using the symmetry, we see that 
\begin{equation}\label{series-heat}
P_tu_0(x) = (P_tu_0,\delta_x)=\sum_y p_t(x,y)u_0(y)
\end{equation}
for $u_0\in\ell^2(V)$. Substituting $u_0=\delta_y$ we see that $p_t(x,y)$ satisfies the heat equation in 
variables $x,t$. We try to extend $P_t$ to functions that are not necessarily in $\ell^2(V)$ by using this
formula and verifying the convergence of the series. To do this we shall need an estimate in the lemma below
of $p_t(x,y)$ for $t\in[0,T]$ and $d(x,y)$ large.
\begin{lemma}\label{heat-decay}
For every $T>0$ there exist a constant $C=C(a,M,T)>0$ such that $$
p_t(x,y) \leq \frac{C}{d(x,y)!}$$ for all $t\in [0,T] $.
\end{lemma}
\begin{proof}
Write $A^n(x,y)$ for the matrix coefficient of the $n$-th
power of $A$. Then $A(x,y)=0$ if $d(x,y)>1$ by the
locality of $A$. It follows, that
$A^n(x,y)=0$ if $d(x,y) > n$. Now suppose that $d(x,y)=k$. It
follows from (\ref{semi}) that
\begin{equation}\label{series}
p_t(x,y) = \sum_{n=k}^\infty\frac{(-t)^nA^n(x,y)}{n!}.
\end{equation}
Since the operator $A$ is bounded with $\parallel A \parallel \leq 2aM$, $$ |A^n(x,y)| =
|(A^n\delta_x,\delta_y)| \leq 2^na^nM^n.$$ Therefore the series
obtained by factoring out $1/k!$ from (\ref{series}) is easily seen
to be uniformly bounded for $t\leq T$. This proves the lemma.
\end{proof}
The lemma says that for $t$ bounded, the heat kernel $p_t(x,y)$ decays very rapidly as the distance
$d(x,y)$ goes to infinity. This is a familiar behavior of the heat kernel of a Riemannian manifold of bounded geometry. Thus we can substitute for $u_0$ in (\ref{series-heat}) functions of moderate growth
so that the series defining $u(x,t)$ converges and produces a solution of (\ref{initial}). In particular,
this yields existence of bounded solutions of (\ref{initial}) asserted in Theorem \ref{unique}. More precisely,
for bounded initial data $|u_0| \leq c$, we define the solution of (\ref{initial}) by (\ref{series-heat}) and group the terms as follows $$
\sum_y p_t(x,y)u_0(y) = \sum_{k=0}^{\infty} \left ( \sum_{d(x,y)=k } p_t(x,y)u_0(y) \right ).$$
By our assumption on the valence, the number of terms in the inner sum is at most $M^k$. Thus, for a bounded $t$, the absolute value of the $k$-th term together with its $t$ derivative is dominated by $(C/k!)M^k c$
because of Lemma \ref{heat-decay}. 
This shows that the series converges very rapidly and can be differentiated term by term proving existence in
Theorem \ref{max-par}.
For future reference we make the following
\begin{remark}\label{allow-growth}
In the argument above we could have allowed $u_0$ to grow  at a certain rate. For example, the argument
goes through if $|u_0(y)| \leq c_1e^{c_2 d(x,y)}$.
\end{remark}

Our next result gives a relation between a ground state and the heat semigroup. It illustrates a
technique used frequently in the study of diffusions \cite{sullivan-lambda}, \cite{pinsky}, 
\cite{dod-mat3}. Let $\calH = \{ u:V \longrightarrow \bbC \mid u\cdot \phi \in \ell^2 (V) \}$.
It is a Hilbert space with the inner product $<u,v>=\sum_{x\in V} u(x)\ol{v}(x)\phi^2(x)$.
We use the ground state $\phi$ to transplant the semigroup $P_t$ to $\calH$. Namely, define $\tilde{P}_t$ 
as a bounded self-adjoint operator on $\calH$ by
\begin{equation}\label{renorm}
\tilde{P}_t = e^{\lambda _0 t} [\phi^{-1}]P_t[\phi] = e^{\lambda _0 t} [\phi^{-1}]e^{-tA}[\phi],
\end{equation}
where $\lambda_0 = \lambda_0(A)$ and $[f]$ denotes the operator of multiplication by a function $f$.
Observe that for $u_0\in \calH$
\begin{equation}\label{renorm-kernel}
\tilde{P}_t u_0(x)=  e^{\lambda_0 t}\sum_y \frac{1}{\phi(x)}p_t(x,y)\phi(y)u_0(y)
\end{equation}
by (\ref{series-heat}).
Clearly, $\tilde{P}_t$, $t\geq 0$ is a semigroup with infinitesimal generator 
$$
-\tilde{A} = -[\phi^{-1}] (A-\lambda_0)[\phi].
$$
The following calculation gives a local formula for $\tilde{A}$.
\begin{eqnarray}
\tilde{A}u(x) & = & \phi^{-1}(x) A(\phi u)(x) -\lambda_0 u(x) \nonumber \\
& = & \phi^{-1}(x) \sum_{y\sim x} a_{x,y} \left (\phi(x)u(x) -\phi(y)u(x)\right ) \nonumber \\
 && + \phi^{-1}(x) \sum_{y\sim x} a_{x,y} \left ( \phi(y)u(x) -\phi(y)u(y)\right )-\lambda_0 u(x) \nonumber\\
& = & \lambda_0 u(x) + \sum_{y\sim x} a_{x,y} \frac{\phi(y)}{\phi(x)}\left (u(x) - u(y)\right ) - \lambda_0 u(x)\nonumber\\
& = & \sum_{y\sim x} a_{x,y} \frac{\phi(y)}{\phi(x)} \left (u(x) - u(y)\right ). \label{generator}
\end{eqnarray}
Note that $\tilde{A}$ is different than the local operators considered until now as its coefficients are
not symmetric in $x,y$. We will consider however the initial value problem analogous to (\ref{initial})
for the operator $\tilde{A}$.
\begin{theorem}
Under the assumptions stated in the beginning of this section, the initial value problem
\begin{equation*}
\begin{split}
\tilde{A}u + \frac{\partial u}{\partial t}  &= 0\\
u(x,0) &= u_0(x)
\end{split}
\end{equation*}
has a unique bounded solution $u(x,t)$ for every bounded function $u_0$.
\end{theorem}
\begin{proof}
The proof is completely analogous to the proof of Theorem \ref{unique}. The uniqueness used only the maximum
principle in Lemma \ref{max-par} which in turn depended only on positivity and \emph{not on symmetry} of the coefficients of
the operator $A$. The proof thus applies equally well to the operator $\tilde{A}$ whose coefficients are positive by (\ref{generator}) since the ground state $\phi$ is positive. Similarly, one proves existence for bounded initial data using the formula (\ref{renorm-kernel}) and applying Remark \ref{allow-growth} together with the estimate of Corollary
\ref{growth-gr-state}.\\
\end{proof}
The following corollary is of independent interest. Its special case was used to derive certain estimates of the
heat kernel for the combinatorial Laplacian in \cite{dod-mat3}.
\begin{corollary}
Under the assumption of this section, the ground state $\phi$ of $A$ is complete i.e.\ satisfies$$
P_t\phi=e^{-\lambda_0t}\phi.$$
\end{corollary}
\begin{proof}
By the theorem above, $\tilde{P}_t$ applied to the function $u_0\equiv 1$ is a solution of the equation
$\tilde{A}u+\frac{\partial u}{\partial t} =0$ with the initial data $u_0$. The function identically equal to
one is also a solution. By uniqueness, the two solutions are equal i.e.
$$
e^{\lambda_0 t}\sum_y \frac{1}{\phi(x)}p_t(x,y)\phi(y)=1
$$ 
for all $t>0,x\in V$.
This proves the corollary.
\end{proof}

{\bfseries Acknowledgement:} I am very grateful to Radek Wojciechowski for a careful reading of the paper, correcting errors and making suggestions that that lead to improvement of exposition.

\end{document}